\documentclass[psamsfont,a4paper,11pt]{amsart}

\usepackage[english]{babel}   

\usepackage{setspace}     
\usepackage{enumerate} 	
\usepackage{afterpage} 	
\usepackage[usenames,dvipsnames]{color}
\usepackage{graphicx}	

\usepackage{amsmath} 	
\usepackage{amssymb}	
\usepackage{mathrsfs} 	
\usepackage{amsfonts} 	
\usepackage{latexsym} 	
\usepackage[latin1]{inputenc} 

\usepackage{amsthm} 	
\usepackage{stackrel}       
\usepackage{amscd} 	

\usepackage{tabularx}	
\usepackage{longtable}	

\usepackage{verbatim} 
\usepackage{blindtext} 

\usepackage{lipsum}

\newcommand\blfootnote[1]{%
  \begingroup
  \renewcommand\thefootnote{}\footnote{#1}%
  \addtocounter{footnote}{-1}%
  \endgroup
}

\allowdisplaybreaks

\newcommand{\C}{\mathbb{C}}

\newcommand{\R}{\mathbb{R}}

\newcommand{\Z}{\mathbb{Z}} 
\newcommand{\Q}{\mathbb{Q}}

\renewcommand{\to}{\longrightarrow}

\def \inf{\mathop{\rm{inf}}}

\newtheorem{Thm}{Theorem}[section]		

\newtheorem{Prop}[Thm]{Proposition}

\theoremstyle{definition}

\newtheorem*{defn}{Definition}  

\theoremstyle{remark}
 
\newtheorem*{rmk}{Remark}

\newtheorem{ind}[]{{\rm\it Indice}}

\usepackage{multicol}

\usepackage{tikz}
\usetikzlibrary{arrows}

\usepackage{paralist}
\usepackage{cite}

\title{Harmonic Maass form eigencurves}

\author[Wagner]{Ian Wagner}

\begin{document}

\maketitle

\blfootnote{$2010$ \textit{Mathematics Subject Classification}.  $11$F$03$, $11$F$37$.}

\begin{abstract}
We construct two families of  harmonic Maass Hecke eigenforms.  Using these families, we construct $p$-adic harmonic Maass forms in the sense of Serre.  The $p$-adic properties of these forms answer a question of Mazur about the existence of an ``eigencurve-type" object in the world of harmonic Maass forms. 

\end{abstract}

\section{Introduction and statement of results}
In \cite{Se1} Serre introduced the notion of a $p$-adic modular form and showed the power of studying a $p$-adic analytic family of modular eigenforms.  Work of Hida \cite{Hida} and Coleman \cite{Cole} expanded on Serre's initial definition of $p$-adic modular form to introduce overconvergent modular forms and offered more examples and applications.  Coleman, in particular, defined the slope of an eigenform as the $p$-adic valuation of its $U_{p}$ eigenvalue and proved that overconvergent modular forms with small slope are classical modular forms.  In \cite{CM} Coleman and Mazur organized all of these results by constructing a geometric object called the ``eigencurve."  The \textit{eigencurve} is a rigid-analytic  curve whose points correspond to normalized finite slope $p$-adic overconvergent modular eigenforms.  

Using Kummer's congruences, Serre was able to give the first examples of $p$-adic modular forms.  Let $v_{p}$ be the $p$-adic valuation on $\Q_{p}$.  If $f = \sum a(n) q^{n} \in \Q[[q]]$ is a formal power series in $q$ then define $v_{p}(f) := \inf_{n}v_{p}(a(n))$.  We then say that $f$ is a \textit{$p$-adic modular form} if there exists a sequence of classical modular forms $f_{i}$ of weights $k_{i}$ such that $v_{p}(f-f_{i}) \to \infty$ as $i \to \infty$.  The weight of a $p$-adic modular form is given by the limits of weights of the classical (holomorphic) modular forms in $X := \Z_{p} \times \Z/(p-1) \Z$.  For a more in-depth discussion of weights see \cite{Se1}.  

The first examples Serre offered came from the Eisenstein series.  Let $\sigma_{k}(n) := \sum_{d \vert n} d^{k}$ be the divisor function, $z = x + iy \in \mathbb{H}$, and $q = e^{2 \pi i z}$.  Then for $k \geq 1$, the weight $2k$ \textit{Eisenstein series} is given by 
\[G_{2k}(z) := \frac{1}{2} \zeta(1-2k) + \sum_{n=1}^{\infty} \sigma_{2k-1}(n) q^{n}, \tag{1.1}\]
where $\zeta(s)$ is the Riemann zeta function.  For $2k \geq 4$, $G_{2k}(z)$ is a weight $2k$ holomorphic modular form on $SL_{2}(\Z)$.  Using the Eisenstein series Serre constructed the $p$-adic Eisenstein series.  Define \[\sigma_{k}^{(p)} := \sum_{\substack{d \vert n \\ \gcd(d, p)=1}} d^{k},\tag{1.2}\]
and let $\zeta^{(p)}(s)$ be the \textit{$p$-adic zeta function} (see \cite{KL}).  We now have that 
\[G_{2k}^{(p)}(z) = \frac{1}{2} \zeta^{(p)}(1-2k) + \sum_{n=1}^{\infty} \sigma_{2k-1}^{(p)}(n) q^{n} \tag{1.3}\]
is a $p$-adic Eisenstein series of weight $2k$.  Clearly there is a sequence $2k_{i}$ of positive even integers that tends to $2k$ $p$-adically and $\sigma_{2k_{i}-1}(n)$ tends to $\sigma_{2k-1}^{(p)}(n)$ $p$-adically.  The $p$-adic Eisenstein series are also classical modular forms on $\Gamma_{0}(p)$ and can be written as
\begin{equation*} 
G_{2k}^{(p)}(z) = G_{2k}(z) - p^{2k-1}G_{2k}(pz).
\end{equation*}
This form is a $p$-stabilization of $G_{2k}(z)$ so that $G_{2k}^{(p)}(z)$ is an eigenform for the $U_{p}$ operator with eigenvalue coprime to $p$.  The $p$-adic Eisenstein series satisfy incredible congruences; we have that $G_{k_{1}}^{(p)}(z) \equiv G_{k_{2}}^{(p)}(z) \pmod{p^{a}}$ whenever $k_{1} \equiv k_{2} \pmod{(p-1)p^{a-1}}$ and $k_{1}$ and $k_{2}$ are not divisible by $p-1$.  For example, $6 \equiv 10 \pmod{4}$ and $6, 10 \not \equiv 0 \pmod{4}$, so we have
\begin{align*}
G_{6}^{(5)}(z) &= \frac{781}{126} + q + 33q^2 + 244q^3 + 1057q^4 +q^5 + \cdots, \end{align*}
and
\begin{align*}
G_{10}^{(5)}(z) &= \frac{488281}{66}+ q + 513q^2 + 19684q^3 + 262657q^4 +q^5+ \cdots \end{align*}
are congruent modulo $5$.  The congruences can be explained using Kummer's congruences and Euler's theorem.  

Mazur recently raised the question of whether or not an eigencurve-like object exists in the world of harmonic Maass forms.  Harmonic Maass forms are traditionally built using methods which rarely lead to forms which are eigenforms (for background see \cite{AMS}).  Namely, the most well known constructions involve Poincar\'{e} series, indefinite theta functions, and Ramanujan's mock theta functions.  These methods do not generally offer Hecke eigenforms.  To this end, the first goal is to construct canonical families of harmonic Hecke eigenforms, out of which one hopes to be able to construct an eigencurve.

Here we construct two families, one integer weight and one half-integer weight, of harmonic Maass forms which are eigenforms for the Hecke operators (see Section $2$ for the definition of the relevant Hecke operators).  We define the weight $k$ differential operator $\xi_{k}$ by
\[\xi_{k} := 2i y^{k} \overline{ \frac{\partial }{\partial \overline{z}}}.\]
The $\xi$-operator defines a surjective map from the space of weight $2-k$ harmonic Maass forms on $\Gamma$ to the space of weight $k$ weakly holomorphic modular forms on $\Gamma$ (see \cite{cbms}).  A natural place to look for a suitable family of harmonic Maass forms is the pullback under the $\xi$-operator of the classical Eisenstein series that Serre used.  The pullback, however, is infinite dimensional.  For example, the $\xi$-operator annihilates weakly holomorphic modular forms.  Therefore, the problem is to construct forms that are the pullback under the $\xi$-operator, and are also Hecke eigenforms and have $p$-adic properties.  Our first family will be a pullback of the classical Eisenstein series that saitsfies these properties.   For $\rm{Re}(s) >0$ or $\rm{Im}(z) >0$, let $\Gamma(s,z) := \int_{z}^{\infty} t^{s-1} e^{-t} dt$ be the \textit{incomplete gamma function}. 
For $k >0$, define \begin{align*}
G(z, -2k) &:= \frac{(2k)! \zeta(2k+1)}{(2 \pi)^{2k}} + \frac{(-1)^{k+1} y^{1+2k} 2^{1+2k} \pi \zeta(-2k-1)}{2k+1} \\
&+ (-1)^{k}(2 \pi)^{-2k} (2k)! \sum_{n=1}^{\infty} \frac{\sigma_{2k+1}(n)}{n^{2k+1}} q^{n} \\
&+ (-1)^{k} (2 \pi)^{-2k} \sum_{n=1}^{\infty} \frac{\sigma_{2k+1}(n)}{n^{2k+1}} \Gamma(1+2k, 4 \pi n y) q^{-n}. \tag{1.4} 
\end{align*}

For half-integral weights, the analogue of the classical Eisenstein series are the Cohen-Eisenstein series \cite{C}.  For more information on half-integral weight modular forms see \cite{cbms}.  Our family of forms will be a pullback of the Cohen-Eisenstein series under the $\xi$-operator.   Define
\[T_{r}^{\chi}(v) := \sum_{a \vert v} \mu(a) \chi(a) a^{r-1} \sigma_{2r-1}(v/a),\]
where $\mu(a)$ is the M\"{o}bius function and $r$ is an integer.  Set $(-1)^{r}N = Dv^2$ with $D$ the discriminant of $\Q(\sqrt{D})$ and let $\chi_{D} = \left( \frac{D}{\cdot} \right)$ be the associated character.  Let
\[c_{r}(N) = \begin{cases} i^{2r+1}L(1+r, \chi_{D})  \frac{1}{v^{2r+1}} T_{r+1}^{\chi_{D}}(v) & N >0 \\
 i^{2r-1} \zeta(1+2r) + \frac{2^{2r+4} i \pi^{2r+1} y^{r+\frac{1}{2}} \zeta(-1-2r)}{(2r-3)\Gamma(2r+1)} & N=0 \\
\pi^{3/2} \frac{L(-r, \chi_{D}) T_{r+1}^{\chi_{D}}(v)}{N^{r + \frac{1}{2}}}  \frac{\Gamma \left( \frac{r+a}{2} \right)}{\Gamma \left(\frac{r + 1+a}{2} \right) \Gamma \left( r +\frac{1}{2} \right)} \Gamma \left(r + \frac{1}{2}, -4 \pi Ny \right) & N<0, 
\end{cases} \tag{1.5}\]
where $a = 0$ if $r$ is odd and $a=1$ if $r$ is even.  Then, for $r \geq 1$, define 
\[\mathcal{H}\left(z, -r+\frac{1}{2}\right) := \sum_{N \in \Z} c_{r}(N) q^N. \tag{1.6}\]
\begin{rmk}
The coefficients for $N>0$ and $N<0$ of $\mathcal{H}\left(z, -r+\frac{1}{2}\right)$ alternate between $L$-functions for real and imaginary quadratic fields as $r$ changes parity.  The $L$-functions for real quadratic fields are known to encode information about the torsion groups of $K$-groups for real quadratic fields.  Therefore, the functions $\mathcal{H}\left(z, -r+\frac{1}{2}\right)$ create a grid that encodes this information for $K_{n}(\Q(\sqrt{D}))$ as both $n$ and $D$ vary.
\end{rmk}

For $k \in \R$, the \textit{weight} $k$ \textit{hyperbolic Laplacian operator} on $\mathbb{H}$ is defined by
\[\Delta_{k} := -y^{2} \left( \frac{\partial^{2}}{\partial x^{2}} + \frac{\partial^{2}}{\partial y^{2}} \right) + iky \left( \frac{\partial}{\partial x} + i \frac{\partial}{\partial y} \right) = -4y^{2} \frac{\partial}{\partial z} \frac{\partial}{\partial \overline{z}} + 2iky \frac{\partial}{\partial \overline{z}}. \tag{1.7}\]
A weight $k$ \textit{harmonic Maass form} on a subgroup $\Gamma$ of $SL_{2}(\Z)$ is a smooth function $f: \mathbb{H} \to \C$ such that it transforms like a modular form of weight $k$ on $\Gamma$, it is annihilated by the weight $k$ hyperbolic Laplacian operator, and it satisfies suitable growth conditions at all cusps.  In particular, we consider harmonic Maass forms with manageable growth, which are defined in Section $2.1$ below.  One subgroup of particular interest is
\[\Gamma_{0}(N):= \Bigg\{ \left( \begin{array}{cc} a & b \\ c & d \end{array} \right) \in SL_{2}(\Z) : c \equiv 0 \pmod{N} \Bigg\}.\]  For a more thorough background on harmonic Maass forms see \cite{AMS}.  We now have the following theorem.

\begin{Thm} \label{Main}
Assuming the notation above, the following are true.
\begin{enumerate}
\item For positive integers $k$, we have that $G(z, -2k)$ is a weight $-2k$ harmonic Maass form on $SL_{2}(\Z)$.  Furthermore, $G(z, -2k)$ has eigenvalue $1 + \frac{1}{p^{2k+1}}$ under the Hecke operator $T(p)$.
\item For positive integers $r$, we have that $\mathcal{H} \left(z, -r +\frac{1}{2} \right)$ is a weight $-r + \frac{1}{2}$ harmonic Maass form on $\Gamma_{0}(4)$.  Furthermore, $\mathcal{H} \left(z, -r +\frac{1}{2} \right)$ has eigenvalue $1 + \frac{1}{p^{2r+1}}$ under the Hecke operator $T(p^{2})$.
\end{enumerate}
\end{Thm}

\begin{rmk}
The proof of Theorem \ref{Main} will show that these forms can be viewed as two parameter functions in $z$ and $w$ where $w$ is the weight of the form.  Specializing $w$ to $-2k$ for the integer weight case and to $-r + \frac{1}{2}$ in the half-integral weight case produces two families of harmonic Maass Hecke eigenforms which define lines on two Hecke eigencurves.
\end{rmk}

\begin{rmk}
Just as the weight $2$ Eisenstein series is not a modular form, the weight $0$ form here is not a harmonic Maass form.  However, we will see that there is a weight $0$ $p$-adic harmonic Maass form in the same way that there is a weight $2$ $p$-adic Eisenstein series.
\end{rmk}

\begin{rmk}
Integer weight non-holomorphic Eisenstein series have been studied before.  For example, in \cite{Za2} Zagier considers the form
\begin{equation*}
\widetilde{G}(z, s) = \frac{1}{2} \sum_{\substack{ (m,n) \in \Z^{2} \\ (m,n) \neq (0,0)}} \frac{y^{s}}{|nz + m|^{2s}},
\end{equation*}
which transforms as a weight $0$ modular form with respect to $z$ is an eigenform of $\Delta_{0}$ with eigenvalue $s(1-s)$.  This form plays an important role in the Rankin-Selberg method \cite{R}, \cite{S}.  Zagier shows that it has a meromorphic continuation so that $\widetilde{G}^{*}(z,s) = \pi^{-s} \Gamma(s) \widetilde{G}(z,s)$ satisfies $\widetilde{G}^{*}(z,s) = \widetilde{G}^{*}(z, 1-s)$.  The Maass lowering operator $L = - 2iy^{2} \frac{\partial}{\partial \overline{z}}$ takes a function that transforms like a modular form of weight $k$ to a function that transforms like a modular form of weight $k-2$.  Furthermore, if $f$ is an eigenform for $\Delta_{k}$ with eigenvalue $\lambda$, then $L(f)$ is an eigenform for $\Delta_{k-2}$ with eigenvalue $\lambda -k +2$ (Chapter $5$ of \cite{AMS}).  In particular, we can see that
\begin{equation*}
L^{k}(\widetilde{G}(z,s)) = \frac{\Gamma(s+k)}{2 \Gamma(s)} \sum_{\substack{ (m,n) \in \Z^{2} \\ (m,n) \neq (0,0)}} \frac{y^{s+k} (nz+m)^{2k}}{|nz+m|^{2s}}
\end{equation*}
is an eigenform of $\Delta_{-2k}$ with eigenvalue $-(s+k)(s-k-1)$.  Evaluating at $s=k+1$ makes the form harmonic and gives the same forms as the ones in Theorem \ref{Main} $(1)$. 

\end{rmk}

\begin{rmk}
The case of $r=0$ has been constructed by Rhoades and Waldherr in \cite{RW} using a slightly different method.  Their result can be recovered using the same method as in this paper and then sieving to suitably modify the Fourier expansion.  The work of Rhoades and Waldherr follows up on work of Duke and Imamo$\overline{\rm{g}}$lu (\cite{DI1}) and Duke, Imamo$\overline{\rm{g}}$lu, and T\'{o}th (\cite{DI2}).  In \cite{DI1} Duke and Imamo$\overline{\rm{g}}$lu use the Kronecker limit formula to construct a function which has values of $L$-functions at $s=1$ for its Fourier coefficients.  This function was the first example and the motivation for the work in \cite{DI2} where Duke, Imamo$\overline{\rm{g}}$lu, and T\'{o}th construct forms of weight $\frac{1}{2}$ on $\Gamma_{0}(4)$ whose Fourier coefficients are given in terms of cycle integrals of the modular $j$-function.
\end{rmk}

\begin{rmk}
The forms in part $1$ of Theorem \ref{Main} behave nicely under the flipping operator (see \cite{AMS}).  Similar functions are studied by Bringmann, Kane, and Rhoades in \cite{BKR}.
\end{rmk}

Serre used the classical Eisenstein series to build $p$-adic modular forms.  In a similar way we can use these harmonic Maass forms to build $p$-adic harmonic Maass forms.
\begin{defn}
A \textbf{weight $k$ $p$-adic harmonic Maass form} is a formal power series  
\[f(z) = \sum_{n \gg -\infty} c_{f}^{+}(n) q^{n} + c_{f}^{-}(0) y^{1-k} + \sum_{0 \neq n \ll \infty} c_{f}^{-}(n) \Gamma \left( 1-k, -4 \pi n y \right) q^{n},\]
where $\Gamma(1-k, -4 \pi ny)$ is taken as a formal symbol and where the coefficients $c_{f}^{\pm}(n)$ are in $\C_{p}$, such that there exists a series of harmonic Maass forms $f_{i}(z)$, of weights $k_{i}$, such that the following properties are satisfied:
\begin{enumerate}
\item $\lim_{i \to \infty} n^{1-k_{i}} c_{f_{i}}^{\pm}(n) = n^{1-k}c_{f}^{\pm}(n)$ for $n \neq 0$.
\item $\lim_{i \to \infty} c_{f_{i}}^{\pm}(0) = c_{f}^{\pm}(0)$.
\end{enumerate}
\end{defn}
\begin{rmk}
Here  $\lim_{i \to \infty} n^{1-k_{i}} c_{f_{i}}^{\pm}(n) = n^{1-k}c_{f}^{\pm}(n)$ means $v_{p}(n^{1-k_{i}} c_{f_{i}}^{\pm}(n) -  n^{1-k}c_{f}^{\pm}(n))$ tends to $\infty$ and we have that $k$ is the limit of the $k_{i}$ in $X$.
\end{rmk}

 We will need a few definitions before describing our $p$-adic harmonic Maass forms.  Let $L_{p}( s, \chi)$ be the \textit{$p$-adic $L$-function} (see \cite{Iw}) and define  
\[T_{r}^{\chi, (p)}(v) := \sum_{\substack{a \vert v \\ \gcd(a, p) = 1}} \mu(a) \chi(a) a^{r-1} \sigma_{2r-1}^{(p)}(v/a). \tag{1.8}\]
Also define the usual \textit{$p$-adic Gamma function} (see \cite{M}) by
\begin{equation*}
\Gamma^{(p)}(n) := (-1)^{n} \prod_{\substack{0<j<n \\ p \nmid j}} j \qquad \text{if } n \in \Z,
\end{equation*}
and  
\begin{equation*}
\Gamma^{(p)}(x) := \lim_{n \to x} \Gamma^{(p)}(n) \qquad \text{if } x \in \Z_{p}.
\end{equation*}
For any $x \in \Z_{p}$ we have $v_{p}(\Gamma^{(p)}(x)) =1$.  In the following formulas we define $\pi :=\Gamma^{(p)} \left( \frac{1}{2} \right)^{2}$ so that $v_{p}(\pi)=1$.  We now have the following theorem.
\begin{Thm} \label{p}
Suppose $p$ is prime and let $\Gamma(\cdot, \cdot)$ be a formal symbol.  Then the following are true.
\begin{enumerate}
\item For each $k \in X$, we have that \begin{align*}
G^{(p)}(z, -2k) &:= \frac{\Gamma^{(p)}(2k+1) \zeta^{(p)}(2k+1)}{(2 \pi)^{2k}} + \frac{(-1)^{k+1} y^{1+2k} 2^{1+2k} \pi \zeta^{(p)}(-2k-1)}{2k+1} \\
&+ (-1)^{k}(2 \pi)^{-2k} \Gamma^{(p)}(2k+1) \sum_{n=1}^{\infty} \frac{\sigma_{2k+1}^{(p)}(n)}{n^{2k+1}} q^{n} \\
&+ (-1)^{k} (2 \pi)^{-2k} \sum_{n=1}^{\infty} \frac{\sigma_{2k+1}^{(p)}(n)}{n^{2k+1}} \Gamma(1+2k, 4 \pi n y) q^{-n} 
\end{align*}
is a weight $-2k$ $p$-adic harmonic Maass form.
\item For each $-r + \frac{1}{2} \in X$, let
\[c_{r}^{(p)}(N) := \begin{cases} i^{2r+1}L_{p}(1+r, \chi_{D})  \frac{1}{v^{2r+1}} T_{r+1}^{\chi_{D}, (p)}(v) & N >0 \\
 i^{2r-1} \zeta^{(p)}(1+2r) + \frac{2^{2r+4} i \pi^{2r+1} y^{r+\frac{1}{2}} \zeta^{(p)}(-1-2r)}{(2r-3)\Gamma^{(p)}(2r+1)}& N=0 \\
\pi^{3/2} \frac{L_{p}(-r, \chi_{D}) T_{r+1}^{\chi_{D}, (p)}(v)}{N^{r + \frac{1}{2}}}  \frac{\Gamma^{(p)} \left( \frac{r+a}{2} \right)}{\Gamma^{(p)} \left(\frac{r + 1+a}{2} \right) \Gamma^{(p)} \left( r +\frac{1}{2} \right)} \Gamma \left(r + \frac{1}{2}, -4 \pi Ny \right) & N<0. 
\end{cases}\]
Then $\mathcal{H}^{(p)}\left(z, -r + \frac{1}{2}\right) = \sum_{N \in \Z} c_{r}^{(p)} q^{N}$ is a weight $-r + \frac{1}{2}$ $p$-adic harmonic Maass form.
\end{enumerate}
\end{Thm}

\begin{rmk}
Note that these forms enjoy congruences similar to the ones for $p$-adic modular forms because of the generalized Bernoulli number congruences.  Congruences for the holomorphic parts rely on the existence of a $p$-adic regulator for $L$-functions.    
\end{rmk}

\begin{rmk}
When $k$ is an integer, the forms $G^{(p)}(z, -2k)$ satisfy
\begin{equation*}
G^{(p)}(z, -2k) = G(z, -2k) - G(pz, -2k),
\end{equation*}
which is the analogue of the equation above that the classical $p$-adic Eisenstein series satisfy.  This implies that $G^{(p)}(z, -2k)$ is a standard harmonic Maass form on $\Gamma_{0}(p)$.  We do not know a similar formula for the half-integral weight forms.   
\end{rmk}

\begin{rmk}
Suppose $p$ is a prime and consider an infinite sequence of even integers which $p$-adically go to zero (i.e. $\{2p^{t} \}_{t=1}^{\infty}$).  By the proof of Theorem \ref{p}, taking the $p$-adic limit of a series of forms with these weights defines a $p$-adic harmonic Maass form of weight $0$.  As noted above, this is the analogue to the quasimodular form $E_{2}$ which is not quite a modular form, but Serre showed leads to a weight $2$ $p$-adic modular form.  In fact, the weight $0$ $p$-adic harmonic Maass form constructed here is the preimage of Serre's weight $2$ $p$-adic Eisenstein series under the $\xi_{0}$-operator.
\end{rmk}

\begin{rmk}
Theorem \ref{p} implies that the Cohen-Eisenstein series are $p$-adic modular forms in the sense of Serre.  This fact was proven by Koblitz in \cite{K}.
\end{rmk}

Not much is known about harmonic Maass Hecke eigenforms except for the forms constructed here.  The fact that Hecke operators increase the order of singularities at cusps poses a major roadblock in the study of harmonic Maass Hecke eigenforms.  The forms constructed here stand out because this issue doesn't arise.  It is an open question of Mazur to describe what the general structure of a ``mock eigencurve" could be.  For example, are there other branches of the mock eigencurve that connect together other harmonic Maass eigenforms?

In Section $2$ we will give background knowledge on harmonic Maass forms and state some results of Zagier vital to the construction of the half-integral weight forms.  Section $3$ will be dedicated to the construction of the forms in Theorem \ref{Main} and proving that they are Hecke eigenforms.  Section $4$ will be used to discuss the $p$-adic properties of the forms in Theorem \ref{p}.  

\section{Background on harmonic Maass forms and results of Zagier}
\subsection{Basics of harmonic Maass forms}
In this section we let $z = x +iy \in \mathbb{H}$, with $x, y \in \R$.  We denote the space of weight $k$ harmonic Maass forms on $\Gamma$ by $H_{k}(\Gamma)$.  The following details on harmonic Maass forms can be found in Chapter $4$ of \cite{AMS}.   If the growth condition mentioned in the definition of harmonic Maass forms given above is given by
\[f(z) = O(e^{\varepsilon y})\] as $y\to \infty$ for some $\varepsilon >0$, then we say that $f$ is a \textit{weight} $k$ \textit{harmonic Maass form of manageable growth on} $\Gamma$ and we denote this space by $H_{k}^{\rm{mg}}(\Gamma)$.

If $f(z)$ is a weight $k$ harmonic Maass form on a congruence subgroup, $\Gamma \subset SL_{2}(\Z)$, it can be naturally decomposed into its \textit{holomorphic part}, $f^{+}(z)$, and its \textit{non-holomorphic part}, $f^{-}(z)$.  The holomorphic part of a harmonic Maass form is often called a \textit{mock modular form}.  The Fourier expansion of $f$ also naturally splits as 
\[ f(z) = \sum_{n \gg -\infty} c_{f}^{+}(n) q^{n} + c_{f}^{-}(0) y^{1-k} + \sum_{\substack{n \ll \infty \\ n \neq 0}}c_{f}^{-}(n) \Gamma(1-k, -4 \pi n y) q^{n}. \tag{2.2}\]
We have the following propostion about the action of the Hecke operators on $f(z)$.
\begin{Prop} [Proposition $7.1$ of \cite{AMS}]
Suppose that $f(z) \in H_{\kappa}^{\rm{mg}}(\Gamma_{0}(N), \chi)$ with $\kappa \in \frac{1}{2} \Z$.  Then the following are true.
\begin{enumerate}
\item For $m \in \mathbb{N}$, we have that $f \vert T(m) \in H_{\kappa}^{\rm{mg}}(\Gamma_{0}(N), \chi)$.
\item if $\kappa \in \Z$, $\epsilon \in \{ \pm \}$, then, unless $n=0$ and $\epsilon = -$,
\[c_{f\vert T(p)}^{\epsilon}(n) = c_{f}^{\epsilon}(pn) + \chi(p) p^{\kappa -1} c_{f}^{\epsilon} \left(\frac{n}{p} \right).\]
Moreover,
\[c_{f\vert T(p)}^{-}(0) = (p^{\kappa -1} + \chi(p)) c_{f}^{-}(0).\]
\item if $\kappa \in \frac{1}{2} \Z \setminus \Z$, then, with $\epsilon \in \{ \pm \}$ ($n \neq 0$ for $\epsilon = -$), we have that
\[c_{f\vert T(p^{2})}^{\epsilon} (n) = c_{f}^{\epsilon}( p^{2}n) + \chi^{*}(p) \left( \frac{n}{p} \right) p^{\kappa - \frac{3}{2}}c_{f}^{\epsilon} (n) + \chi^{*}(p^2) p^{2 \kappa -2} c_{f}^{\epsilon} \left( \frac{n}{p^2} \right),\]
where $\chi^{*}(n) := \left( \frac{(-1)^{\kappa - \frac{1}{2}}}{n} \right) \chi(n)$.  If $n = 0$ and $\epsilon = -$, then we have that
\[c_{f \vert T(p^2)}^{-} (0) = (p^{-2+2 \kappa} + \chi^{*}(p^2))c_{f}^{-}(0).\]
\end{enumerate}
\end{Prop}

Differential operators are an important tool for studying harmonic Maass forms.  We will focus on the $\xi$-operator.  Let $M_{k}^{!}(\Gamma_{0}(N))$ be the space of weakly holomorphic modular forms on $\Gamma_{0}(N)$ (see \cite{cbms}).
\begin{Prop} [Theorem $5.10$ of \cite{AMS}]
For any $k \geq 2$, we have that
\[\xi_{2-k}: H_{2-k}^{\rm{mg}}(\Gamma_{0}(N)) \twoheadrightarrow M_{k}^{!}(\Gamma_{0}(N)).\]
In particular, for $f \in H_{2-k}^{\rm{mg}}(\Gamma_{0}(N))$, we have that 
\[\xi_{2-k}(f(z)) = \xi_{2-k}(f^{-}(z)) = (k-1)\overline{c_{f}^{-}(0)}-(4 \pi)^{k-1} \sum_{n \gg - \infty} \overline{c_{f}^{-}(-n)} n^{k-1} q^{n}.\]
\end{Prop}
The $\xi$-operator allows for a connection between the Hecke operators for harmonic Maass forms and modular forms (see \cite{AMS}).  In particular, we have that
\[p^{d(1- \kappa)} \xi_{\kappa}(f \vert T(p^{d}, \kappa, \chi)) = \xi_{\kappa}(f) \vert T(p^{d}, 2-\kappa, \chi), \tag{2.3}\]
where 
\[ d := \begin{cases} 1 & \rm{if} \ \kappa \in \Z, \\ 2 & \rm{if} \ \kappa \in \frac{1}{2} + \Z. \end{cases}\]

\subsection{Results of Zagier}
Several results of Zagier will be applicable to the construction of our forms.  We will state them here.

\begin{Prop}[Zagier, \cite{Za1}] \label{Dirichlet}
For positive integers $a$ and $c$, let 
\[\lambda(a,c) = \begin{cases} i^{\frac{1-c}{2}} \left( \frac{a}{c} \right) & \emph{if } c \emph{ is odd, } a \emph{ even} \\ 
i^{\frac{a}{2}} \left( \frac{c}{a} \right) & \emph{if } a \emph{ is odd, } c \emph{ even} \\
0 & \rm{otherwise.} \end{cases} \tag{2.4}\]
Define the Gauss sum $\gamma_{c}(n)$ by 
\[ \gamma_{c}(n):=\frac{1}{\sqrt{c}} \sum_{a=1}^{2c} \lambda(a,c) e^{- \pi i n \frac{a}{c}}. \tag{2.5}\]
Let $n$ be a nonzero integer and define a Dirichlet series $E_{n}(s)$ by
\[E_{n}(s):= \frac{1}{2} \sum_{\substack{c=1 \\ c \ odd}}^{\infty} \frac{\gamma_{c}(n)}{c^{s}} + \frac{1}{2} \sum_{\substack{c=2 \\ c \ even}}^{\infty} \frac{\gamma_{c}(n)}{(c/2)^{s}}, \tag{2.6}\]
(i.e. $E_{n}(s) = \sum a_{m} m^{-s}$ where $a_{m} = \frac{1}{2} (\gamma_{m}(n) + \gamma_{2m}(n))$ when $m$ is odd, and $a_{m} = \frac{1}{2} \gamma_{2m}(n)$ when $m$ is even).
Let $K = \Q(\sqrt{n})$, $D$ be the discriminant of $K$, $\chi_{D} = \left(\frac{D}{\cdot} \right)$ be the character of $K$, and $L(s, \chi_{D}) = \sum \frac{\chi_{D}(n)}{n^{s}}$ be the $L$-series of $K$ (if $n$ is a perfect square, then $\chi(m) = 1$ for any $m$ and $L(s, \chi) = \zeta(s)$).  Then if $n \equiv 2, 3 \pmod{4}$, we have 
\[E_{n}(s) = 0.\]
If $n \equiv 0, 1 \pmod{4}$, we have 
\[E_{n}(s) = \frac{L(s, \chi_{D})}{\zeta(2s)} \sum_{\substack{a, c \geq 1 \\ ac \vert v}} \frac{\mu(a) \chi_{D}(a)}{c^{2s-1} a^{s}} = \frac{L(s, \chi_{D})}{\zeta(2s)} \frac{T_{s}^{\chi_{D}}(v)}{v^{2s-1}}, \tag{2.7}\]
where $n = v^{2}D$ and 
\[T_{s}^{\chi}(v) = \sum_{t \vert v} t^{2s-1} \sum_{a \vert t} \frac{\mu(a) \chi(a)}{a^s} = \sum_{a \vert v} \mu(a) \chi(a) a^{s-1} \sigma_{2s-1}(v/a). \tag{2.8}\]
Furthermore, we have 
\[E_{0}(s) = \frac{\zeta(2s-1)}{\zeta(2s)}. \tag{2.9}\]

\end{Prop}

\begin{rmk}
It is clear from Zagier's proof in \cite{Za1} that $E_{n}(s)$ can be continued to a meromorphic function on the whole $s$-plane.  It will also be beneficial to note that $T_{s}^{\chi}(v) = v^{2s-1} T_{1-s}^{\chi}(v)$.  
\end{rmk}
It will be useful for us to define 
\[E_{n}^{odd}(s) = \sum_{\substack{c = 1 \\ c \ odd}}^{\infty} \gamma_{c}(n) c^{-s}, \tag{2.10}\]
and 
\[E_{n}^{even}(s) = \sum_{\substack{c=1 \\ c \ even}}^{\infty} \gamma_{c}(n) (c/2)^{-s}, \tag{2.11}\]
so that \[E_{n}(s) = \frac{1}{2} \left( E_{n}^{odd}(s) + E_{n}^{even}(s) \right). \tag{2.12}\]

\section{Proof of Theorem \ref{Main}}
Here we prove Theorem \ref{Main}.  There are two cases to consider, the integer weight and half-integral weight cases.  In the next subsection we consider the integer weight case.
\subsection{Proof of Theorem \ref{Main} Part $1$}
We will construct the forms from Theorem \ref{Main} part $1$ first.  Let $z \in \mathbb{H}$ and $k \in \Z$.  Define
\[\mathcal{G}(z, -2k, s) := \frac{1}{2} \sum_{n,m}{}^{\prime} \frac{(mz + n)^{2k}}{\vert mz+n \vert^{2s}}, \tag{3.1} \]
where the primed sum means the sum runs over all $(n,m)$ except $(0,0)$.  $\mathcal{G}(z, -2k, s)$ has a meromorphic continuation to the whole $s$-plane.  Let
\[f(z, -2k, s) := \sum_{n = - \infty}^\infty (z+n)^{2k} \vert z+n \vert^{-2s}.  \tag{3.2}\]
Then we have \[f(z, -2k, s) = \sum_{n=-\infty}^\infty h_{n,2k}(s,y) e^{2 \pi i n z} = \sum_{n=-\infty}^\infty h_{n}(y, -2k, s) e^{2 \pi i n x} e^{-2 \pi n y},\]
where \[h_{n}(y, -2k, s) = \int_{iy-\infty}^{iy + \infty} z^{2k} \vert z \vert^{-2s} e^{-2 \pi i n z} dz.\]
After making the substitution $z = yt + iy$ we have
\[h_{n}(y, -2k, s) = y^{1 +2k -2s} e^{2 \pi n y} \int_{-\infty}^\infty (t + i)^{2k} (t^2 + 1)^{-s} e^{-2 \pi i ny t} dt.  \tag{3.3}\]
For $n=0$, we have $h_{0}(y, -2k, s) = y^{1+2k-2s} \int_{-\infty}^{\infty} (t+i)^{2k}(t^{2}+1)^{-s} dt$.  Following Zagier, we choose our branch cut along the negative imaginary axis.  Then using contour integration, we find that
\[h_{0}(y, -2k, s) = 2i y^{1+2k-2s} e^{k \pi i} \sin(\pi (s-2k)) \int_{-i}^{-i \infty} \vert t+i \vert^{2k-s} \vert t-i \vert^{-s} dt.  \tag{3.4}\]
We substitute $t$ for $-i(2u +1)$ to arrive at 
\[h_{0}(y, -2k, s) = 2^{2+2k-2s} y^{1+2k-2s}e^{k \pi i} \sin(\pi (s-2k)) \int_{0}^{\infty} u^{2k-s} (u+1)^{-s} du.\]
We make one more substitution, $u = \frac{1-v}{v}$.  Then, we have
\begin{align*} 
h_{0}(y, -2k, s)&= 2^{2+2k-2s} y^{1+2k-2s} e^{k \pi i} \sin(\pi (s-2k)) \int_{0}^{1} (1-v)^{2k-s} v^{2s-2k-2} dv \\
&=  2^{2+2k-2s} y^{1+2k-2s} e^{k \pi i} \pi \frac{\Gamma(2s-2k-1)}{\Gamma(s-2k) \Gamma(s)}.  \tag{3.5}
\end{align*}
For $n>0$, we define a path $c_{1}$ as a clockwise path around $-i$ from $-i \infty$ to $-i \infty$.  Then we have
\[h_{n}(y, -2k, s) = y^{1+2k-2s} \int_{c_{1}} (v+i)^{2k} (v^2+1)^{-s} e^{-2 \pi i n y v} dv.\]
Substitute $v$ for $t-i$ and define the path $c_{2} = c_{1} + i$, then we have 
\[h_{n}(y, -2k, s) = y^{1+2k-2s} e^{-2 \pi n y} \int_{c_{2}} t^{2k-s} (t-2i)^{-s} e^{-2 \pi i n y t} dt.\]
For $n<0$, define the path $c_{3}$ as before to circle $i$ clockwise from $i \infty$ to $i \infty$.  Making the substitutions $v = t+i$ and $c_{4} = c_{3} - i$, we arrive at 
\[h_{n}(y, -2k, s) = y^{1+2k-2s} e^{2 \pi n y} \int_{c_{4}} t^{-s} (t+2i)^{2k-s} e^{-2 \pi i n y t} dt, \]
for $n < 0$.  Notice that $h_{n}(my, -2k, s) = m^{1+2k-2s} h_{mn}(y, -2k, s)$, so we have
\begin{align*}
\mathcal{G}(z, -2k, s) &= \frac{1}{2} \sum_{m=- \infty}^\infty f(mz, -2k, s) \\
&= \zeta(2s-2k) + \sum_{m=1}^\infty f(mz, -2k, s) \\
&= \zeta(2s-2k) + \sum_{m=1}^\infty \sum_{n=- \infty}^\infty m^{1+2k-2s} h_{mn}(y, -2k, s) e^{2 \pi i nm x}. \tag{3.6}
\end{align*}
We want to now look at the limit as $s$ goes to zero in order to obtain a negative weight Eisenstein series.  However, it is clear that for any $n \in \Z$, $h_{n}(y, -2k, 0) = 0.$  Thus, our $\mathcal{G}(z, -2k, 0)$ functions will also go to zero.  In order to work around this we will look at the derivative of our Eisenstein series with respect to $s$.  Define
\[G(z, -2k) := \lim_{s \to 0} \frac{d}{ds} \mathcal{G}(z, -2k, s). \tag{3.7}\]
We will now calculate the $q$-expansion of $G(z, -2k)$.  For $n=0$, we have 
\begin{align*}
\frac{d}{ds} h_{0}(y, -2k, s)\vert_{s=0} &= y^{1+2k} 2^{2+2k} e^{k \pi i} \pi \frac{\Gamma(-2k-1)}{\Gamma(-2k)} \\
&= \frac{(-1)^{k+1} y^{1+2k} 2^{2+2k} \pi}{2k+1}. \tag{3.8}
\end{align*}
For $n>0$, we have
\begin{align*}
\frac{d}{ds} h_{n}(y, -2k, s)\vert_{s=0} &= -y^{1+2k} e^{-2 \pi n y} \int_{c_{2}} t^{2k} \log(t(t-2i)) e^{-2 \pi i nyt} dt \\
&= -y^{1+2k} e^{-2 \pi ny} (2 \pi i) \int_{-i \infty}^{0} t^{2k} e^{-2 \pi i nyt} dt \\
&= (-1)^{k} y^{1+2k} e^{-2 \pi ny}(2 \pi ) \int_{0}^\infty t^{2k} e^{-2 \pi n yt} dt \\
&= (-1)^{k}(2 \pi)^{-2k} n^{-2k-1} \Gamma(2k+1) e^{-2 \pi ny}.  \tag{3.9}
\end{align*}
The log term jumps by $2 \pi i$ across the branch cut, while everything else is continuous.  Similarly, for $n<0$ we have
\[\frac{d}{ds} h_{n}(y, -2k, s)\vert_{s=0} = (-1)^{k+1} (2 \pi)^{-2k} n^{-2k-1} e^{-2 \pi ny} \Gamma(1+2k, -4 \pi n y), \tag{3.10}\]
Let $h_{n}'(y, -2k, 0) := \frac{d}{ds} h_{n}(y, -2k, s)\vert_{s=0}$, then, from equation $2.6$, we have 
\[G(z, -2k) = 2 \zeta'(-2k) + \sum_{n=- \infty}^\infty h_{n}'(y, -2k, 0) \sigma_{2k+1}(n) e^{2 \pi i n x}.  \tag{3.11}\]
Recall that $\sigma_{2k+1}(0) = \frac{1}{2} \zeta(-2k-1)$.  Putting everything together leads to the construction of the forms in Theorem \ref{Main} part $1$.  A short calculation shows these forms are harmonic.  In order to show that $G(z, -2k)$ is a Hecke eigenform, notice that its image under the $\xi$-operator is a nonzero multiple of the weight $2k+2$ Eisenstein series, $E_{2k+2}(z)$.  $E_{2k+2}$ is known to be an eigenform with eigenvalue $\sigma_{2k+1}(p) = 1 + p^{2k+1}$ under the Hecke operator $T(p)$.  By equation ($2.3$) and inspection it is clear that $G(z, -2k)$ is then an eigenform with eigenvalue $1 + \frac{1}{p^{2k+1}}$.
{\flushright \qed}

\subsection{Proof of Theorem \ref{Main} part $2$}
Let $k = 2r-1$ with $r \geq 1$.  We define the two Eisenstein series $F \left(z, -\frac{k}{2}, s \right)$ and $E \left(z, -\frac{k}{2}, s \right)$ by
\[F \left(z, -\frac{k}{2}, s\right) = \sum_{\substack{n,m \in \Z \\ n > 0 \\ 4 \vert m}} \left(\frac{m}{n} \right) \varepsilon_{n}^{-k} \frac{(mz +n)^{k/2}}{\vert mz+n \vert^{2s}}, \tag{3.12}\]
and
\[E\left(z, -\frac{k}{2}, s \right) = \frac{(2z)^{k/2}}{\vert 2z \vert^{2s}} F\left(\frac{-1}{4z}, -\frac{k}{2}, s \right),\]
where $\left(\frac{m}{n} \right)$ is the \textit{Kronecker symbol} and 
\[\varepsilon_{n} := \begin{cases} 1 & if \ n \equiv 1 \pmod{4} \\ i & if \ n \equiv 3 \pmod{4}. \end{cases}\]   A linear combination of these forms will have a meromorphic continuation to the whole $s$-plane and evaluating at $s=0$ will give our weight $-\frac{k}{2}$ form.  We will abuse this fact by letting $s=0$ in the assembly of the forms.  We have 
\[E\left(z, -\frac{k}{2}, s\right) = 2^{\frac{k}{2} - 2s} \sum_{\substack{ n,m \in \Z \\ n>0, odd}} \left(\frac{m}{n} \right) \varepsilon_{n}^{-k} \frac{(nz -m)^{k/2}}{\vert nz-m \vert^{2s}}. \tag{3.13}\]
From this we have 
\begin{align*}
E \left(z, -\frac{k}{2}, s \right)  &= 2^{\frac{k}{2} - 2s} \sum_{n>0, odd} \varepsilon_{n}^{-k} n^{\frac{k}{2} -2s} \sum_{m \pmod{n}} \left( \frac{m}{n} \right) \sum_{h = - \infty}^{\infty} \frac{ \left(z - \frac{m}{n} + h \right)^{\frac{k}{2}}}{\vert z - \frac{m}{n} + h \vert^{2s}} \\
&= 2^{\frac{k}{2} - 2s} \sum_{n>0, odd} \varepsilon_{n}^{-k} n^{\frac{k}{2} -2s} \sum_{N=- \infty}^{\infty} \sum_{m \pmod{n}} \left( \frac{m}{n} \right) \alpha_{N}\left(y, -\frac{k}{2}, s\right) e^{- \frac{2 \pi i N m}{n}} q^{N} \\
&= \sum_{N= - \infty}^{\infty} a(N) q^{N}, 
\end{align*}
where 
\[a(N) = 2^{\frac{k}{2} - 2s} \alpha_{N}\left(y, -\frac{k}{2}, s \right) \sum_{n>0, odd} \varepsilon_{n}^{-k} n^{\frac{k}{2} -2s} \sum_{m \pmod{n}} \left( \frac{m}{n} \right) e^{- \frac{2 \pi i N m}{n}}, \tag{3.14}\]
and by the Poisson summation formula
\[\alpha_{N}\left(y, -\frac{k}{2}, s\right) = \int_{iy - \infty}^{iy + \infty} z^{\frac{k}{2}} \vert z \vert^{-2s} e^{-2 \pi i N z} dz.\]
Making the substitution $z = yt + iy$ gives us
\[\alpha_{N}\left(y, -\frac{k}{2}, s\right) = y ^{\frac{k}{2} +1 -2s} e^{2 \pi N y} \int_{- \infty}^{\infty} (t + i)^{\frac{k}{2}} (t^{2} + 1)^{-s} e^{-2 \pi i N y t} dt.\]
Following Zagier, we choose the branch cut along the negative imaginary axis.  Using contour integration we have
\[\alpha_{N}\left(y, -\frac{k}{2}, s \right) = 2 e^{\frac{k \pi i}{4}} \sin \left( \pi \left( \frac{k}{2} - s \right) \right) y^{\frac{k}{2} + 1 - 2s} \int_{0}^{\infty} t^{\frac{k}{2} - s} (t+2)^{-s} e^{-2 \pi N y t} dt. \tag{3.16}\]
Letting $s=0$ we arrive at \begin{align*}
\alpha_{N}\left(y, -\frac{k}{2}, s\right) &= 2 e^{\frac{k \pi i}{4}} \sin \left( \frac{\pi k}{2} \right) y^{\frac{k}{2} + 1} \int_{0}^{\infty} t^{\frac{k}{2}} e^{-2 \pi N y t} dt \\
&= 2 e^{\frac{k \pi i}{4}} \sin \left( \frac{\pi k}{2} \right) y^{\frac{k}{2} + 1} (2 \pi N y)^{-\frac{k}{2} - 1} \int_{0}^{\infty} t^{\frac{k}{2}} e^{-t} dt \\
&= 2 e^{\frac{k \pi i}{4}} \sin \left( \frac{\pi k}{2} \right) (2 \pi N)^{- \frac{k}{2} - 1} \Gamma \left(\frac{k}{2} + 1 \right), \tag{3.17}
\end{align*}
for $N >0$.  If we evaluate the similar integral for $N \leq 0$, because we do not cross a branch cut the integral is zero.  It will be useful to evaluate the derivative.  We have that
\begin{align*}
\frac{d}{ds} \alpha_{N}\left(y, -\frac{k}{2}, s\right) \vert_{s=0} &= - y^{\frac{k}{2} + 1} e^{4 \pi N y} (2 \pi i) \int_{i \infty}^{0} (t + 2i)^{\frac{k}{2}} e^{- 2 \pi i N y t} dt \\
&= -2 y^{\frac{k}{2}+1}\pi i^{\frac{k}{2}} \int_{2}^{\infty} t^{\frac{k}{2}} e^{2 \pi Nyt} dt \\
&= i^{-\frac{k}{2}} (2 \pi)^{- \frac{k}{2}} N^{-\frac{k}{2}-1} \Gamma \left( \frac{k}{2} + 1, -4 \pi Ny \right), \tag{3.18}
\end{align*}
for $N < 0$, while 
\begin{equation*} \frac{d}{ds} \alpha_{0}\left(y, -\frac{k}{2}, s\right) \vert_{s=0} = -\frac{2^{\frac{7}{2} -r} i^{\frac{k}{2}}y^{\frac{k}{2}+1} \pi}{2r-3}. \end{equation*}  Similarly, for $F\left(z, -\frac{k}{2}, s \right)$ we have
\begin{align*}
F\left(z, -\frac{k}{2}, s \right) &= 1+ \sum_{\substack{m>0 \\ 4 \vert m}} m^{\frac{k}{2} - 2s} \sum_{n \pmod{m}} \left( \frac{m}{n} \right) \varepsilon_{n}^{-k} \sum_{h= -  \infty}^{\infty} \frac{\left( z + \frac{n}{m} + h \right)^{\frac{k}{2}}}{ \vert z + \frac{n}{m} + h \vert^{2s}} \\
&= 1+ \sum_{\substack{m>0 \\ 4 \vert m}} m^{\frac{k}{2} - 2s} \sum_{n \pmod{m}} \left( \frac{m}{n} \right) \varepsilon_{n}^{-k} \sum_{N = -\infty}^{\infty} \alpha_{N}\left(y, -\frac{k}{2}, s \right) e^{\frac{2 \pi i N n}{m}}q^{N} \\
&= 1+ \sum_{N= - \infty}^{\infty} b(N) q^{N},
\end{align*}
where 
\[b(N) = \alpha_{N}\left(y, -\frac{k}{2}, s \right) \sum_{\substack{m>0 \\ 4 \vert m}} m^{\frac{k}{2} - 2s} \sum_{n \pmod{m}} \left( \frac{m}{n} \right) \varepsilon_{n}^{-k} e^{\frac{2 \pi i N n}{m}}. \tag{3.19}\]
Using Proposition \ref{Dirichlet} and by manipulating the inner sums of $a(N)$ and $b(N)$, it is not hard to show that
\begin{align*} a(N) &= 2^{\frac{k}{2} - 2s} \alpha_{N}\left(y, -\frac{k}{2}, s \right) \sum_{n>0, odd} n^{\frac{k}{2} -2s} \sum_{\substack{m = 1 \\ m \ even}}^{2n} \lambda(m, n) e^{- \pi i (-1)^{r} N \frac{m}{n}} \\
&= 2^{\frac{k}{2} - 2s} \alpha_{N}\left(y, -\frac{k}{2}, s \right) \sum_{n>0, odd} n^{\frac{k}{2} + \frac{1}{2} -2s} \gamma_{n}((-1)^{r}N) \\
&= 2^{\frac{k}{2} + 1 - 2s} \alpha_{N}\left(y, -\frac{k}{2}, s \right) \frac{1}{2} E_{(-1)^{r}N}^{odd} \left(-\frac{k}{2} - \frac{1}{2} + 2s \right), \tag{3.20}
\end{align*}
and 
\begin{align*}
b(N) &= (1 + i^{2r+1}) 4^{\frac{k}{2} + \frac{1}{2} -2s} \alpha_{N}\left(y, -\frac{k}{2}, s\right) \frac{1}{2} \sum_{\substack{m>0 \\ m \ even}} \frac{\gamma_{m}((-1)^{r}N)}{(m/2)^{-\frac{k}{2} - \frac{1}{2} + 2s}} \\
&= (1 + i^{2r+1}) 4^{\frac{k}{2} + \frac{1}{2} -2s} \alpha_{N}\left(y, -\frac{k}{2}, s \right) \frac{1}{2} E_{(-1)^{r}N}^{even} \left(-\frac{k}{2} - \frac{1}{2} + 2s \right). \tag{3.21}
\end{align*}

We are now able to define our forms.  Define
\begin{align*} \mathcal{H}\left(z, -r +\frac{1}{2} \right)& = \sum_{N=- \infty}^{\infty} c_{r}(N) q^{N} \\
&:= \lim_{s \to 0} \zeta(1 + 2r-4s) \left( i^{2r-1}F\left(z, -r + \frac{1}{2}, s \right) + 2^{r - \frac{1}{2}}(1 + i^{2r-1}) E\left(z, -r + \frac{1}{2}, s \right) \right). \tag{3.22}\end{align*}
The rest of the construction is using Proposition \ref{Dirichlet}.  Similar calculations can be found in \cite{C} or \cite{Za2}.  Note that the functional equations for the zeta function and the $L$-function are used and that there is  pole when evaluating the non-holomorphic coefficients.  The image of $\mathcal{H} \left(z, -r +\frac{1}{2} \right)$ under the $\xi$-operator is a nonzero multiple of the weight $r+\frac{3}{2}$ Cohen-Eisenstein series.  The weight $r + \frac{3}{2}$ Cohen-Eisenstein series is a Hecke eigenform with eigenvalue $1 + p^{2r+1}$ under the Hecke operator $T(p^{2})$.  Therefore, using equation ($2.3$), we can see that
\[\mathcal{H} \left(z, -r + \frac{1}{2} \right) \Big |T(p^{2}) - \left(1 + \frac{1}{p^{2r+1}} \right) \mathcal{H} \left(z, -r +\frac{1}{2} \right)\]
is a weight $-r +\frac{1}{2}$ holomorphic modular form in the Kohnen plus space (see \cite{cbms}).  This space is empty and so $\mathcal{H} \left(z, -r +\frac{1}{2} \right)$ must be a Hecke eigenform with eigenvalue $1 + \frac{1}{p^{2r+1}}$.
{\flushright \qed}

\section{Proof of Theorem \ref{p}}
In order to discuss $p$-adic harmonic Maass forms we will first need to recall some facts about Bernoulli numbers.  Values of the Reimann zeta function at negative integers are tied to Bernoulli numbers.  In fact we have $\zeta(1-2k) = -\frac{B_{2k}}{2k}$, and $\zeta(-2k)=0$.  In a similar way, there is a connection between generalized Bernoulli numbers and the values of $L$-functions at negative integers.  The \textit{generalized Bernoulli numbers} $B(n, \chi)$ are defined by the generating function
\[\sum_{n=0}^\infty B(n, \chi) \frac{t^{n}}{n!} = \sum_{a=1}^{m-1} \frac{\chi(a) t e^{at}}{e^{mt}-1}, \tag{4.1} \]
Where $\chi$ is a Dirichlet character modulo $m$.  Generalized Bernoulli numbers are known to give the values of Dirichlet L-functions at non-positive integers.  In fact, from \cite{cbms} we know that if $k$ is a positive integer and $\chi$ is a nontrivial Dirichlet character, then 
\[L(1-k, \chi) = - \frac{B(k, \chi)}{k}.\tag{4.2} \]
This connection helps one define a $p$-adic $L$-function, $L_{p}(s, \chi)$.  The $p$-adic $L$-function is analytic except for a pole at $s=1$ with residue $\left(1 - \frac{1}{p} \right)$.  For $n \geq1$ we have that
\[L_{p}(1-n, \chi) = - (1 - \chi \cdot \omega^{-n}(p)p^{n-1}) \frac{B(n, \chi \cdot \omega^{-n})}{n}, \tag{4.3}\]
where $\omega$ is the Teichm\"{u}ller character.  The Teichm\"{u}ller character is a $p$-adic Dirichlet character of conductor $p$ if $p$ is odd and conductor $4$ if $p=2$.  It is best to view it as a $p$-adic object.  For more information see Chapter $5$ of \cite{W}.  Kummer famously showed that if $n \equiv m \pmod{(p-1)p^{a}}$ and $ (p-1) \nmid n,m$ for an odd prime $p$, then 
\[(1 - p^{n-1}) \frac{B_{n}}{n} \equiv (1-p^{m-1}) \frac{B_{m}}{m} \pmod{p^{a+1}}, \tag{4.4}\]
where $a$ is a nonnegative integer.  Similar congruences hold for generalized Bernoulli numbers as well.  For example, if we let $\chi \neq 1$ be a primitive Dirichlet character with conductor not divisible by $p$, then if $n \equiv m \pmod{p^{a}}$ we have
\[(1 - \chi \cdot \omega^{-n}(p)p^{n-1}) \frac{B(n, \chi \cdot \omega^{-n})}{n} \equiv (1- \chi \cdot \omega^{-m}(p)p^{m-1}) \frac{B(m, \chi \cdot \omega^{-m})}{m} \pmod{p^{a+1}}. \tag{4.5}\]
Notice that twisting by the appropriate power of the Teichm\"{u}ller character removes the dependence on the residue class of $n$ and $m$ modulo $p-1$ here.  The family of $p$-adic harmonic Maass forms coming from the integer weight forms in Theorem \ref{Main} are constructed  in the exact same way as the $p$-adic Eisenstein series in \cite{Se1}.  Equation $(4.4)$ shows that the constant term, the $p$-adic zeta function at a negative integer, will satisfy congruences.  The other terms satisfy congruences due to Euler's theorem which generalizes Fermat's Little Theorem.  The algebraic parts of these $p$-adic harmonic Maass forms enjoy similar congruences as their modular counterparts.  In fact, the non-holomorphic parts are nearly identical to the $p$-adic Eisenstein series.  The holomorphic parts behave not quite as nicely only because the $p$-adic zeta function at positive integers does not behave as nicely as at negative integers.  However, it is still expected that it satisfies similar congruences modulo some $p$-adic regulator. For example, we have
\begin{align*}&G^{+,(5)}(z , -2) =  -\frac{1}{2 \pi^{2}}\left(\zeta^{(5)}(3) + q + \frac{9}{8}q^{2} + \frac{28}{27} q^{3} + \frac{73}{64} q^{4} + \frac{1}{75} q^{5} + \cdot \cdot \cdot \right), \end{align*}
while 
\[G^{+,(5)}(z, -6) = - \frac{45}{4 \pi^{6}} \left(\zeta^{(5)}(7) + q + \frac{129}{128} q^2 + \frac{2188}{2187} q^3 + \frac{16513}{16384} q^4 + \frac{1}{78125} q^5 + \cdot \cdot \cdot \right).\]
The family of $p$-adic harmonic Maass forms coming from the half-integral weight forms from Theorem \ref{Main} are defined using $p$-adic $L$-functions and the fact that $T_{r}^{\chi, (p)}(v)$ is the $p$-adic limit of $T_{r}^{\chi}(v)$.  As in the previous case, the non-holomorphic parts satisfy nice congruences due to equation $(4.5)$.  The holomorphic parts are expected to satisfy congruences modulo a $p$-adic regulator.

\section*{Acknowledgements}
The author would like to thank Barry Mazur, J-P. Serre, Larry Rolen, Michael Griffin, Kathrin Bringmann, and \"{O}zlem Imamo$\bar{\rm{g}}$lu for their comments on an earlier version of this paper.  The author would also like to thank Ken Ono for his numerous suggestions and the two referees for their comments which improved the quality of this paper.

\end{document}